\documentclass[]{interact}

\usepackage{epstopdf}
\usepackage[caption=false]{subfig}

\usepackage[numbers,sort&compress]{natbib}
\bibpunct[, ]{[}{]}{,}{n}{,}{,}
\renewcommand\bibfont{\fontsize{10}{12}\selectfont}
\usepackage{amsmath,amssymb}
\usepackage{hyperref}
\usepackage{geometry}
\usepackage{enumitem}
\usepackage{enumitem}
\usepackage{algorithm}
\usepackage{algorithmic}
\usepackage{amsthm}
\usepackage{tikz}
\usepackage{tcolorbox}
\usepackage{framed}
\usepackage{mathtools}

\newtheorem{theorem}{Theorem}[section]
\newtheorem{corollary}[theorem]{Corollary}
\newtheorem{proposition}[theorem]{Proposition}
\theoremstyle{definition}
\newtheorem{definition}[theorem]{Definition}
\theoremstyle{remark}
\newtheorem{remark}[theorem]{Remark}

\usepackage{amsthm}

\newtheorem{lemma}{Lemma}[section]

\newtheorem*{lemmarep}{Lemma A.1 (Common–limit constraint)}
\newtheorem*{propsync}{Proposition (Averaging synchronization)}

\begin{document}

\articletype{RESEARCH ARTICLE}

\title{Piecewise Recursive Sequences with Adaptive Thresholds : Boundary Convergence and Applications}

\title{Piecewise Recursive Sequences with Adaptive Thresholds : Boundary Convergence and Applications}

\author{
\name{Slimane Alaoui Soulimani Valenti\textsuperscript{a}\thanks{CONTACT Slimane Alaoui Soulimani Valenti. Email: slimane.alaoui@mail.utoronto.ca}}
\affil{\textsuperscript{a}University of Toronto, Toronto, Canada}
}

\maketitle

\begin{abstract}
We study discrete-time dynamical systems that switch between different evolution rules based on thresholds that themselves adapt over time. Specifically, we analyze the coupled recursion $a_{n+1} = f(a_n)$ if $a_n \leq c_n$ and $a_{n+1} = g(a_n)$ if $a_n > c_n$, where the threshold evolves according to $c_{n+1} = h(a_n, c_n)$. By transforming the system into triangular coordinates, we map the problem to a piecewise-smooth system with a fixed switching boundary. We derive explicit local stability criteria based on the lower-triangular structure of the associated Jacobians and establish a ``common-limit constraint'': we prove that any convergent orbit that switches regimes infinitely often must converge to a limit where the state and threshold coincide. To demonstrate the framework's utility, we develop an asset-pricing model where investor sentiment follows a ``trailing-stop'' rule. We characterize the parameter regions for market stability versus bubble formation and numerically extend the analysis to coupled markets, illustrating how adaptive thresholds can facilitate the propagation of instability and contagion across financial networks.
\end{abstract}

\section{Introduction}

Many discrete-time systems switch between regimes according to a threshold that itself adapts over time. We study the coupled recursion
\[
a_{n+1}=
\begin{cases}
f(a_n), & a_n \le c_n,\\[2pt]
g(a_n), & a_n > c_n,
\end{cases}
\qquad
c_{n+1}=h(a_n,c_n),
\]
in which both the state \(a_n\) and the switching threshold \(c_n\) evolve endogenously. This formulation appears in applications in which decision rules react to outcomes while the decision criteria are themselves updated.

A key structural observation is that the change of variables
\[
x_n := a_n, \qquad y_n := a_n - c_n
\]
reduces the model to a \emph{triangular}, piecewise-defined map with a fixed switching boundary \(y=0\):
\[
x_{n+1}=
\begin{cases}
f(x_n), & y_n \le 0,\\
g(x_n), & y_n > 0,
\end{cases}
\qquad
y_{n+1}=
\begin{cases}
f(x_n)-h(x_n,x_n-y_n), & y_n \le 0,\\
g(x_n)-h(x_n,x_n-y_n), & y_n > 0.
\end{cases}
\]
Thus the “dynamic threshold’’ is equivalently encoded by a fixed partition with discontinuity along a line, and the \(x\)-component does not depend on \(y\). This places the problem within established lines of work on discontinuous/triangular maps and border effects, including foundational treatments of piecewise-smooth dynamics and border-collision phenomena \cite{dibernardo2008,Gardini_2016}. \cite{nusse1992}

\paragraph*{Position in the literature.}
There is by now a substantial literature on piecewise-smooth and border-collision maps, with general overviews and many examples in \cite{dibernardo2008,Gardini_2016,simpson2016}.  Previous contributions already showed how simple one-dimensional recursions with switching or discontinuity can display a broad range of behaviours \cite{nusse1992,Keener1980,elaydi2005}, and more recent work has developed detailed bifurcation and stability analyses, including higher-dimensional and application-driven models \cite{avrutin2010,avrutin2019,dieci2010,GARDINI2022483}.  Related perspectives come from the study of piecewise monotone and nonlinear difference equations and their invariant measures or periodicities \cite{lasota1973,grove2004,sedaghat2003}, from monotone trajectories in multivalued dynamical systems \cite{aubin1977}, and from the broader theory of nonsmooth and piecewise-linear systems \cite{leine2006,pettit1997,kunze1997,colombo2011,novaes2015,glendinning2017}.  In parallel, the hybrid- and switched-systems literature provides complementary tools for stability and Lyapunov analysis \cite{branicky1998,johansson1998,goebel2012,liberzon2003}, and work on coupled map lattices offers a discrete-time viewpoint on spatially extended interactions \cite{kaneko1992}.

Our use of a triangular planar map is close in spirit to this fixed-boundary and hybrid-systems literature: the switching condition is still given by a scalar inequality, but we keep the switching level as an explicit state variable instead of treating it as a parameter.  In that sense, the present paper can be read as a modest complement to existing results on piecewise-smooth maps and related nonsmooth systems \cite{dibernardo2008,Gardini_2016,simpson2016,lasota1973,grove2004,sedaghat2003,leine2006,branicky1998,johansson1998,goebel2012,liberzon2003,colombo2011,dercole2005,novaes2015,glendinning2017}.

On the applications side, our example in Section~\ref{sec:economic} follows the heterogeneous-agent asset-pricing tradition in which traders switch between simple forecasting or trading rules and the resulting feedback generates nonlinear price dynamics.  Classical references include \cite{Chiarella1992,Lux1995,BrockHommes1998}, and there is a broad macro- and finance-oriented literature on nonlinear and threshold dynamics in economic models \cite{DayHuang1990,hansen2000,hommes2013}.  More recent analyses of models with switching, discontinuities and border collisions in financial contexts can be found in \cite{DieciWesterhoff2013,DieciHe2018,GARDINI2022483,GARDINI2025104436}.  The sentiment-threshold mechanism considered here is a stylized variant of this type of switching rule, with the difference that the threshold is updated endogenously.  Our aim is not to propose a new empirical specification, but rather to use a very simple asset-pricing example to illustrate how the qualitative properties of an adaptive threshold can be studied within a planar piecewise-smooth and hybrid-systems framework.

\paragraph*{Scope and contributions.}
The scope of the paper is fairly modest.  We focus on a simple class of coupled recursions in which a state variable and a switching threshold evolve jointly according to a triangular update rule, and we collect some basic qualitative results that seem useful for this type of model.  More concretely:
\begin{enumerate}
  \item We study a planar map in which the state and the threshold are updated together.  Within this setting we isolate a basic ``common-limit'' observation for bounded orbits that cross the switching boundary infinitely often: any such orbit must have limit points on the boundary where the state and the threshold coincide.  This elementary fact is used later to prove convergence statements and to relate the original coupled formulation to a reduced (fixed-boundary) formulation, in the spirit of the standard theory for piecewise-smooth maps and related nonsmooth systems \cite{dibernardo2008,Gardini_2016,Keener1980,avrutin2010,avrutin2019,elaydi2005,simpson2016,GARDINI2022483,lasota1973,grove2004,sedaghat2003,leine2006,colombo2011,novaes2015,glendinning2017}.  The viewpoint is also compatible with seeing the model as a very simple hybrid or switched system \cite{branicky1998,goebel2012,liberzon2003,johansson1998}.
  \item For the triangular map, we derive one-sided linear conditions for local stability of both interior and boundary fixed points, expressed in terms of the eigenvalues of the branch-wise Jacobians.  These conditions are specialized to our set-up but consistent with the usual stability and border-collision considerations in the piecewise-smooth and non-smooth dynamics literature \cite{dibernardo2008,Gardini_2016,simpson2016,dieci2010,dercole2005,leine2006,kunze1997,pettit1997}.  The results are intended as simple tools for our examples rather than as a general classification.
  \item For a tractable class of adaptive rules (in particular, convex averaging of the current state and the past threshold), we show that whenever the state converges to a fixed point $a^\ast$, the adaptive threshold converges to the same limit, so that $|a_n-c_n|\to 0$.  This synchronization property indicates that, in these cases, the adaptive threshold does not create additional long-run equilibria beyond those of the underlying one-dimensional recursion.  It sits comfortably alongside existing work on asymptotic behaviour, invariant sets and periodicity in difference equations and related discrete-time systems \cite{lasota1973,grove2004,sedaghat2003,aubin1977,kaneko1992}.
  \item Finally, we embed the framework into a stylized heterogeneous-agent asset-pricing model with a dynamic sentiment threshold, in the tradition of \cite{Chiarella1992,Lux1995,BrockHommes1998,DayHuang1990,hommes2013,DieciWesterhoff2013,DieciHe2018,GARDINI2022483,GARDINI2025104436}.  The example is deliberately simple and is used only to illustrate how the general observations on convergence, local stability and synchronization can help organize the analysis of models where switching rules are driven by an evolving threshold; we do not attempt any calibration or detailed empirical evaluation.
\end{enumerate}

\section{Preliminaries and Notation}\label{sec:prelim}

\paragraph*{Original formulation.}
We consider the coupled recursion on an interval $I\subseteq\mathbb{R}$:
\begin{equation}\label{eq:original}
a_{n+1} =
\begin{cases}
f(a_n), & \text{if } a_n \le c_n,\\[2pt]
g(a_n), & \text{if } a_n > c_n,
\end{cases}
\qquad
c_{n+1} = h(a_n,c_n),
\end{equation}
where $f,g:I\to I$ and $h:I\times I\to I$.
Unless explicitly stated otherwise, all objects are real-valued and scalar. Throughout, indices $n\in\mathbb{N}$ start at $n=0$.

\paragraph*{Triangular reduction.}
Introduce the change of variables
\begin{equation}\label{eq:change}
x_n := a_n, \qquad y_n := a_n - c_n .
\end{equation}
Then the switching condition $a_n\le c_n$ becomes $y_n\le 0$, and \eqref{eq:original} is equivalent to the piecewise-defined, triangular map
\begin{equation}\label{eq:triangular}
(x_{n+1},y_{n+1}) =
\begin{cases}
\big(f(x_n),\; f(x_n)-h(x_n,x_n-y_n)\big), & \text{if } y_n \le 0,\\[3pt]
\big(g(x_n),\; g(x_n)-h(x_n,x_n-y_n)\big), & \text{if } y_n > 0,
\end{cases}
\end{equation}
on the state space $X:=I\times \mathbb{R}$. By construction, the $x$-component does not depend on $y$ (triangular structure), while the switching boundary is the fixed line
\[
\Gamma := \{(x,y)\in X : y=0\}.
\]
We will refer to the two half-planes
\[
X^- := \{(x,y)\in X : y\le 0\}, \qquad
X^+ := \{(x,y)\in X : y> 0\}.
\]
The two branches of \eqref{eq:triangular} will be denoted by $T^-:X^-\to X$ and $T^+:X^+\to X$, with the full map $T:X\to X$ defined by $T(z) = T^\pm(z)$ if $z\in X^\pm$, respectively. The reduction \eqref{eq:change}--\eqref{eq:triangular} places the problem within the standard framework of piecewise-smooth maps with a fixed switching manifold; see the general notation used in studies of border-collision and piecewise-smooth dynamics \cite{dibernardo2008,Gardini_2016}.%

\subsection{Standing assumptions}\label{subsec:assumptions}
We work under the following regularity and domain assumptions, which are minimal for the statements made later and align with the fixed-boundary formulation:
\begin{itemize}
  \item[\textbf{(A1)}] $f,g:I\to I$ and $h:I\times I\to I$ are continuous. On any compact subset $K\subset I$ (and $K\times K$ for $h$), the restrictions of $f,g,h$ are Lipschitz.
  \item[\textbf{(A2)}] Away from the boundary $\Gamma$, each branch $T^\pm$ of \eqref{eq:triangular} is $C^1$ on its domain of definition. Derivatives at points on $\Gamma$ are understood in the appropriate one-sided sense when needed.
  \item[\textbf{(A3)}] Unless stated otherwise, all orbits considered are well-defined for all $n\ge 0$ and remain in $X$ (forward invariance of the relevant domain).
\end{itemize}
Assumption \textbf{(A2)} ensures that local linearization, spectral radius, and hyperbolicity statements are meaningful on each smooth branch, in the sense commonly adopted for piecewise-smooth systems \cite{dibernardo2008}.%

\paragraph{Notation}
For a function $h : \mathbb{R}^2 \to \mathbb{R}$ we denote by
\[
\partial_1 h(a,c) := \frac{\partial h}{\partial a}(a,c)
\quad\text{and}\quad
\partial_2 h(a,c) := \frac{\partial h}{\partial c}(a,c)
\]
the partial derivatives of $h$ with respect to its first and second argument, respectively.
In particular, when we write $\partial_2 h(a^*,a^*)$ we mean the derivative of $h$ with respect to its second variable, evaluated at the point $(a^*,a^*)$.

\subsection{Orbits, switching times, and limit sets}\label{subsec:orbits}
Given an initial condition $z_0=(x_0,y_0)\in X$, the forward orbit is $\{z_n\}_{n\ge 0}$ with $z_{n+1}=T(z_n)$. We record the \emph{switching-time sets}
\[
\mathcal{N}^- := \{\,n\ge 0 : z_n\in X^- \,\}, 
\qquad
\mathcal{N}^+ := \{\,n\ge 0 : z_n\in X^+ \,\}.
\]
We say the orbit \emph{visits both partitions infinitely often} if both $\mathcal{N}^-$ and $\mathcal{N}^+$ are infinite. Convergence, boundedness, and $\omega$-limit sets are used in the usual sense. In particular, we will be concerned with cases where $\{x_n\}$ and/or $\{y_n\}$ converge, and with consequences of infinite switching under convergence of the relevant components.

\subsection{Fixed points and local terminology}\label{subsec:fixedpoints}
A point $z^\ast=(x^\ast,y^\ast)\in X$ is a \emph{fixed point} of $T$ if $T(z^\ast)=z^\ast$. We distinguish:
\begin{itemize}
  \item \emph{Interior fixed points in a branch:} $z^\ast\in \operatorname{int}(X^\pm)$ with $T^\pm(z^\ast)=z^\ast$.
  \item \emph{Boundary fixed points:} $z^\ast\in \Gamma$ with the one-sided condition $T^\pm(z^\ast)=z^\ast$ for at least one admissible side (the side depends on whether $y^\ast=0$ is approached from $X^-$ or $X^+$).
\end{itemize}
When $z^\ast\in \operatorname{int}(X^\pm)$ and \textbf{(A2)} holds, the Jacobian $DT^\pm(z^\ast)$ is well-defined and, by the triangular structure of \eqref{eq:triangular}, has the form
\[
DT^\pm(z^\ast) \;=\;
\begin{pmatrix}
\partial_x \Phi^\pm(x^\ast) & 0\\[2pt]
\ast & \partial_y \Psi^\pm(x^\ast,y^\ast)
\end{pmatrix},
\]
where $\Phi^-(x)=f(x)$, $\Phi^+(x)=g(x)$, and
\[
\Psi^\pm(x,y) =
\begin{cases}
f(x)-h(x,x-y), & (\pm = -),\\
g(x)-h(x,x-y), & (\pm = +).
\end{cases}
\]
In particular, $DT^\pm(z^\ast)$ is lower triangular, so its eigenvalues are the diagonal entries. We say that an interior fixed point is \emph{hyperbolic} if the spectral radius of $DT^\pm(z^\ast)$ is different from $1$. For boundary fixed points, we use the standard one-sided linearization terminology from piecewise-smooth dynamics, explicitly indicating the side to which the derivative pertains \cite{dibernardo2008}.%

\subsection{Back-translation to the original variables}\label{subsec:back}
The correspondence \eqref{eq:change} is bijective between $(a_n,c_n)$ and $(x_n,y_n)$. For any statement expressed in $(x,y)$, the back-translation is
\[
a_n = x_n, \qquad c_n = x_n - y_n .
\]
In particular, the boundary $y=0$ corresponds to the diagonal $\{(a,c)\in I\times I: a=c\}$ in the original variables. When convenient, we will express limits and fixed points in whichever coordinates make the statement more transparent.

\medskip
\noindent\textbf{Remark.}
Section~\ref{sec:prelim} only fixes notation and terminology. All substantive statements (e.g., conditions for local stability on a side, or consequences of infinite switching under convergence) are given precisely in the subsequent sections. Our use of one-sided derivatives and hyperbolicity follows standard practice for border-collision and piecewise-smooth maps \cite{dibernardo2008}; see also related applications and bifurcation discussions in \cite{Gardini_2016}, and economic switching contexts in \cite{GARDINI2022483}.

\section{Fixed Points and Local Stability}\label{sec:fixed}

\subsection{Characterization of fixed points}\label{subsec:char-fp}
Let $T:X\to X$ be the triangular map in \eqref{eq:triangular} with branches $T^\pm$ on $X^\pm$.
\begin{definition}[Interior and boundary fixed points]
A point $z^\ast=(x^\ast,y^\ast)\in X$ is a fixed point of $T$ if $T(z^\ast)=z^\ast$.
\begin{itemize}
  \item $z^\ast$ is an \emph{interior fixed point} of the branch $T^\pm$ if $z^\ast\in \operatorname{int}(X^\pm)$ and $T^\pm(z^\ast)=z^\ast$.
  \item $z^\ast$ is a \emph{boundary fixed point} if $z^\ast\in \Gamma:=\{y=0\}$ and $T^\pm(z^\ast)=z^\ast$ for at least one admissible side.
\end{itemize}
\end{definition}

For interior fixed points, using $\Phi^-(x)=f(x)$, $\Phi^+(x)=g(x)$ and
\[
\Psi^\pm(x,y)=
\begin{cases}
f(x)-h(x,x-y),& (\pm=-),\\
g(x)-h(x,x-y),& (\pm=+),
\end{cases}
\]
the fixed point equations are
\begin{equation}\label{eq:fp-interior}
x^\ast=\Phi^\pm(x^\ast),\qquad
y^\ast=\Psi^\pm(x^\ast,y^\ast).
\end{equation}
For boundary fixed points $(x^\ast,0)\in\Gamma$, the one-sided condition $T^\pm(x^\ast,0)=(x^\ast,0)$ reduces to
\begin{equation}\label{eq:fp-boundary}
x^\ast=\Phi^\pm(x^\ast),\qquad h(x^\ast,x^\ast)=x^\ast.
\end{equation}
Thus any boundary fixed point must satisfy $h(x^\ast,x^\ast)=x^\ast$, together with $x^\ast$ being a fixed point of the active branch on the relevant side.

\subsection{Jacobian structure and hyperbolicity}\label{subsec:jac}
On $\operatorname{int}(X^\pm)$, assumption \textbf{(A2)} ensures differentiability of $T^\pm$. By the triangular form of \eqref{eq:triangular}, the Jacobian at $z^\ast\in\operatorname{int}(X^\pm)$ is lower triangular:
\begin{equation}\label{eq:jac}
DT^\pm(z^\ast)=
\begin{pmatrix}
\Phi^{\pm\,\prime}(x^\ast) & 0\\[2pt]
\partial_x \Psi^\pm(x^\ast,y^\ast) & \partial_y \Psi^\pm(x^\ast,y^\ast)
\end{pmatrix}.
\end{equation}
Hence the eigenvalues are the diagonal entries. A direct computation gives, for both branches,
\begin{equation}\label{eq:psy-y}
\partial_y \Psi^\pm(x^\ast,y^\ast)=\partial_2 h\big(x^\ast,\,x^\ast-y^\ast\big),
\end{equation}
while $\Phi^{-\prime}(x^\ast)=f'(x^\ast)$ and $\Phi^{+\prime}(x^\ast)=g'(x^\ast)$. In particular, on $\Gamma$ we have $\partial_y \Psi^\pm(x^\ast,0)=\partial_2 h(x^\ast,x^\ast)$.

\begin{definition}[Hyperbolicity]
An interior fixed point $z^\ast\in\operatorname{int}(X^\pm)$ is \emph{hyperbolic} if $\rho\!\left(DT^\pm(z^\ast)\right)\neq 1$, where $\rho(\cdot)$ denotes the spectral radius. A boundary fixed point is \emph{one-sided hyperbolic} (from side $\pm$) if the one-sided Jacobian along that side has spectral radius different from $1$; see standard usage in piecewise-smooth dynamics \cite{dibernardo2008}.

See also \cite{simpson2016} for border-collision normal form and coexistence phenomena.
\end{definition}

\subsection{Local stability criteria (branch-wise)}\label{subsec:local-stab}
The following branch-wise criteria are immediate consequences of \eqref{eq:jac}--\eqref{eq:psy-y} and the continuity hypotheses; they are stated in a form consistent with one-sided analysis near $\Gamma$ \cite{dibernardo2008}.

\begin{proposition}[Interior, single-branch]\label{prop:interior-stab}
Let $z^\ast\in\operatorname{int}(X^\pm)$ satisfy \eqref{eq:fp-interior} and \textbf{(A2)}. Then the eigenvalues of $DT^\pm(z^\ast)$ are
\[
\lambda_1=\Phi^{\pm\,\prime}(x^\ast),\qquad
\lambda_2=\partial_2 h\big(x^\ast,\,x^\ast-y^\ast\big).
\]
In particular:
\begin{enumerate}
\item If $\max\{|\lambda_1|,|\lambda_2|\}<1$, then $z^\ast$ is (branch-wise) locally asymptotically stable.
\item If $\max\{|\lambda_1|,|\lambda_2|\}>1$, then $z^\ast$ is unstable.
\item The case $\max\{|\lambda_1|,|\lambda_2|\}=1$ is nonhyperbolic and requires separate analysis.
\end{enumerate}
\end{proposition}

\begin{proposition}[Boundary, one-sided]\label{prop:boundary-stab}
Let $(x^\ast,0)\in\Gamma$ satisfy the boundary fixed point conditions \eqref{eq:fp-boundary} for side $\pm$. Suppose $T^\pm$ is $C^1$ on that side near $(x^\ast,0)$. Then the one-sided eigenvalues are
\[
\lambda_1=\Phi^{\pm\,\prime}(x^\ast),\qquad
\lambda_2=\partial_2 h(x^\ast,x^\ast).
\]
If $\max\{|\lambda_1|,|\lambda_2|\}<1$ then $(x^\ast,0)$ is locally asymptotically stable from side $\pm$; if $\max\{|\lambda_1|,|\lambda_2|\}>1$ it is unstable from that side. Nonhyperbolic cases require separate treatment \cite{dibernardo2008}.
\end{proposition}

\begin{remark}[Back-translation]
In the original variables $(a,c)$, a boundary fixed point corresponds to $(a^\ast,c^\ast)$ with $a^\ast=c^\ast$ and $h(a^\ast,a^\ast)=a^\ast$. The branch-wise eigenvalues specialize to $\{f'(a^\ast),\partial_2 h(a^\ast,a^\ast)\}$ or $\{g'(a^\ast),\partial_2 h(a^\ast,a^\ast)\}$, respectively.
\end{remark}

\medskip

\section{A Common-Limit Lemma for Infinite Switching}\label{sec:common-limit}

In this section we isolate a simple but useful constraint on orbits
of the triangular map~\eqref{eq:triangular} that converge while
switching infinitely often across the fixed boundary
$\Gamma=\{y=0\}$. The result is purely kinematic: it does not assert
that such convergence occurs, but shows that any putative limit point
must be a common boundary fixed point for the two branches $T^\pm$.

\begin{definition}[Infinite switching]
Let $\{z_n\}_{n\ge 0}$ be an orbit of $T$, with $z_n=(x_n,y_n)\in X$.
We say that $\{z_n\}$ \emph{visits both partitions infinitely often}
(or \emph{exhibits infinite switching across $\Gamma$}) if the index
sets
\[
\mathcal{N}^- \;:=\; \{n\ge 0 : y_n \le 0\}, \qquad
\mathcal{N}^+ \;:=\; \{n\ge 0 : y_n > 0\}
\]
are both infinite.
\end{definition}

\begin{lemma}[Common-limit constraint]\label{lem:common-limit}
Assume \textbf{(A1)}--\textbf{(A3)}. Let $\{z_n\}$ be an orbit of $T$
that visits both partitions infinitely often and converges to a point
$z^\ast=(x^\ast,0)\in\Gamma$, i.e.\ $z_n\to z^\ast$ as $n\to\infty$.
Then $z^\ast$ is a boundary fixed point for \emph{both} one-sided
branches:
\[
T^-(z^\ast) \;=\; z^\ast
\qquad\text{and}\qquad
T^+(z^\ast) \;=\; z^\ast.
\]
Equivalently, $z^\ast$ satisfies the boundary fixed point conditions
\eqref{eq:fp-boundary} simultaneously for $T^-$ and $T^+$.
\end{lemma}

\begin{proof}
Since $\{z_n\}$ visits both partitions infinitely often, there exist
strictly increasing sequences of indices
\[
\{n_k\}_{k\ge 1} \subset \mathcal{N}^-,
\qquad
\{m_k\}_{k\ge 1} \subset \mathcal{N}^+
\]
with $n_k\to\infty$ and $m_k\to\infty$ as $k\to\infty$.
By assumption $z_n\to z^\ast$, so in particular
\[
z_{n_k} \to z^\ast,
\qquad
z_{m_k} \to z^\ast
\qquad\text{as }k\to\infty.
\]

For $n_k\in\mathcal{N}^-$ we have $z_{n_k}\in X^-$ and therefore
\[
z_{n_k+1} \;=\; T^-(z_{n_k}).
\]
Assumption \textbf{(A2)} ensures that $T^-$ is continuous on $X^-$ up
to the boundary $\Gamma$, so by continuity,
\[
z^\ast \;=\; \lim_{k\to\infty} z_{n_k+1}
\;=\; \lim_{k\to\infty} T^-(z_{n_k})
\;=\; T^-(z^\ast).
\]
The same argument applies to the subsequence $\{m_k\}\subset\mathcal{N}^+$:
for each $m_k$,
\[
z_{m_k+1} \;=\; T^+(z_{m_k}),
\]
and continuity of $T^+$ on $X^+$ up to $\Gamma$ implies
\[
z^\ast \;=\; \lim_{k\to\infty} z_{m_k+1}
\;=\; \lim_{k\to\infty} T^+(z_{m_k})
\;=\; T^+(z^\ast).
\]
Thus $T^-(z^\ast)=T^+(z^\ast)=z^\ast$, so $z^\ast$ is a common
boundary fixed point of the two branches, as claimed.
\end{proof}

The lemma provides only a \emph{necessary} condition: it does not
guarantee that an orbit with infinite switching converges, nor does it
give conditions for global convergence or synchronization. Rather, it
shows that if such convergence occurs, the limit point must be a
boundary point fixed by both $T^-$ and $T^+$.

In the triangular setting of~\eqref{eq:triangular}, the one-sided
branches have the form
\[
T^\pm(x,y) \;=\; \bigl(\Phi^\pm(x),\,\Psi^\pm(x,y)\bigr),
\]
with $\Phi^-(x)=f(x)$, $\Phi^+(x)=g(x)$ and
\[
\Psi^\pm(x,y) =
\begin{cases}
f(x)-h(x,x-y), & (\pm = -),\\[2pt]
g(x)-h(x,x-y), & (\pm = +),
\end{cases}
\]
so the boundary fixed point conditions~\eqref{eq:fp-boundary} become
particularly explicit.

\begin{corollary}[Common limit in triangular coordinates]
Under the hypotheses of Lemma~\ref{lem:common-limit}, suppose
moreover that the orbit is generated by the triangular
system~\eqref{eq:triangular} and converges to $z^\ast=(x^\ast,0)$.
Then the following equalities necessarily hold:
\[
h(x^\ast,x^\ast) \;=\; x^\ast,
\qquad
f(x^\ast) \;=\; x^\ast,
\qquad
g(x^\ast) \;=\; x^\ast.
\]
Equivalently, in the original variables $(a_n,c_n)$ one has
\[
\lim_{n\to\infty} a_n \;=\; \lim_{n\to\infty} c_n,
\]
so any convergent orbit with infinite switching across the fixed
boundary must asymptotically collapse to a single common value.
\end{corollary}

\begin{remark}[Interpretation and relation to the literature]
The fixed-boundary viewpoint isolates the role of the discontinuity
surface from other mechanisms that can generate complicated dynamics.
Lemma~\ref{lem:common-limit} shows that, in the presence of infinite
switching, convergence forces compatibility between the one-sided
dynamics on $X^-$ and $X^+$ at the limiting boundary point. This type
of constraint is implicit in the analysis of many piecewise-smooth
systems (see, e.g., \cite{dibernardo2008} for background on
border-collision bifurcations), and it can be combined with more
detailed model-specific information to restrict possible limit sets or
to interpret numerical observations of apparent synchronization.
\end{remark}

\section{Economic Application: Dynamic Sentiment Thresholds in Financial Markets}\label{sec:economic}

In this section we illustrate how the abstract framework of Section~\ref{sec:prelim}
can be interpreted as a simple asset–pricing model with adaptive
sentiment thresholds. Discrete-time models of financial markets with
heterogeneous agents and nonlinear price dynamics are widely used to
study volatility clustering, bubbles and crashes; see, for example,
\cite{Chiarella1992,Lux1995,BrockHommes1998,DieciWesterhoff2013,DieciHe2018,GARDINI2022483}.
Our goal here is not to propose a fully realistic market microstructure,
but rather to show that a reduced-form “bubble monitor’’ with an adaptive
danger threshold fits exactly into the two-dimensional setting
analyzed in Sections~\ref{sec:fixed}--\ref{sec:common-limit}.

\subsection{Single Market with Trailing-Stop Sentiment Thresholds}\label{subsec:single}

\subsubsection{Mispricing and Adaptive Threshold}

Consider a risky asset with log price $p_t$ at date $t$ and constant
log fundamental value $F$. As in \cite{DieciHe2018,GARDINI2022483}, we work with the
mispricing variable
\begin{equation}\label{eq:mispricing}
  x_t = p_t - F,
\end{equation}
so that $x_t > 0$ corresponds to overvaluation and $x_t < 0$ to
undervaluation.

Investors monitor the level of mispricing relative to an adaptive
“danger’’ or “sentiment’’ threshold $\tau_t$. This threshold represents
the level of mispricing the market has become accustomed to: when
$x_t$ exceeds $\tau_t$, agents perceive the situation as an exuberant
bubble, whereas for $x_t \leq \tau_t$ the market is viewed as roughly
normal. We model $\tau_t$ as a trailing average of past mispricing:
\begin{equation}\label{eq:threshold}
  \tau_{t+1} = (1-\gamma)\tau_t + \gamma x_t, \qquad \gamma\in(0,1),
\end{equation}
so that $\gamma$ controls the speed at which market participants
update their perception of what counts as “normal’’ mispricing.
This has the interpretation of a slowly moving reference point or
“normalization of deviance’’: persistent overvaluation gradually
lifts the danger threshold, while persistent cheapness lowers it.

In the notation of Section~\ref{sec:prelim} this
corresponds to taking $a_t = x_t$, $c_t = \tau_t$ and
\begin{equation*}
  h(x,\tau) = (1-\gamma)\tau + \gamma x,
\end{equation*}
so that $(x_t,\tau_t)$ evolves according to a two-dimensional
triangular map of the type studied in Sections~\ref{sec:fixed}
and~\ref{sec:common-limit}.

\subsubsection{Two Behavioral Regimes: Fundamentalists vs.\ Exuberance}

To keep the dynamics strictly two-dimensional and aligned with the
piecewise-scalar recursion of Section~\ref{sec:prelim}, we work with a
reduced-form price update that depends only on current mispricing and
the current threshold. When the mispricing is below the danger level
($x_t \leq \tau_t$), fundamentals dominate and the market exhibits
mean reversion. When the mispricing exceeds the threshold ($x_t >
\tau_t$), trend-following and exuberant sentiment dominate and push
prices further away from fundamentals. This gives a simple
“trailing-stop’’ or “bubble monitor’’ rule:
\begin{equation}\label{eq:dynamics}
  x_{t+1} =
  \begin{cases}
    \alpha x_t, & \text{if } x_t \leq \tau_t
    \quad\text{(normal / mean-reverting regime)},\\[0.3em]
    \beta x_t + \delta, & \text{if } x_t > \tau_t
    \quad\text{(exuberant / bubble regime)},
  \end{cases}
\end{equation}
with parameters
\[
  0 < \alpha < 1, \qquad \beta > 1, \qquad \delta \in \mathbb{R}.
\]
Here:
\begin{itemize}
  \item $\alpha$ captures the net strength of stabilizing forces
  (fundamentalist trading, liquidity provision) in the normal regime,
  similar in spirit to the stabilizing components in
  \cite{Chiarella1992,Lux1995,BrockHommes1998}.
  \item $\beta$ captures the net strength of destabilizing, trend-following
  or sentiment-driven trading once the bubble is perceived as underway,
  in line with the positive-feedback channels emphasized in
  \cite{BrockHommes1998,DieciHe2018,GARDINI2022483}.
  \item $\delta$ is a constant optimism bias or shift in demand
  when the exuberant regime is active.
\end{itemize}
We emphasize that~\eqref{eq:dynamics} is a reduced-form
law of motion that compresses the detailed trader-level structure
into a piecewise-affine map; this is standard in the heterogeneous
agent literature, where equilibrium price dynamics often reduce to
low-dimensional difference equations
\cite{BrockHommes1998,DieciWesterhoff2013,GARDINI2022483}.

Combining~\eqref{eq:threshold} and~\eqref{eq:dynamics}, the state of
the system at date $t$ is
\[
  z_t := (x_t,\tau_t) \in \mathbb{R}^2.
\]
The switching boundary is given by the line
\[
  \Gamma := \{(x,\tau)\in\mathbb{R}^2 : x = \tau\},
\]
and the evolution is governed by a piecewise-smooth map $T$ that
coincides with a different affine map on each side of~$\Gamma$.

\subsubsection{Identification with the Abstract Model}

Define the scalar maps
\[
  f(x) = \alpha x, \qquad g(x) = \beta x + \delta,
\]
and the threshold update $h$ as above. Then the coupled recursion
\begin{equation*}
  x_{t+1} =
  \begin{cases}
    f(x_t), & x_t \leq \tau_t,\\
    g(x_t), & x_t > \tau_t,
  \end{cases}
  \qquad
  \tau_{t+1} = h(x_t,\tau_t)
\end{equation*}
is exactly of the form studied in Section~\ref{sec:prelim}:
a scalar piecewise recursion for $x_t$ with an
adaptive threshold $\tau_t$ that depends on the current state
$(x_t,\tau_t)$.

The transformation $y_t = x_t - \tau_t$ used in Subsection~\ref{subsec:back}
yields a triangular representation of the dynamics in local
coordinates $(x_t,y_t)$ with switching manifold $y_t = 0$. This is
precisely the setting in which the interior and boundary stability
results of Propositions~\ref{prop:interior-stab}--\ref{prop:boundary-stab}
and the common-limit constraint of Lemma~\ref{lem:common-limit}
apply.

\subsubsection{Jacobian Matrices and Stability}

On each side of the switching line $\Gamma$ the map is affine and
hence has a constant Jacobian. Writing $z_t=(x_t,\tau_t)$ and
$z_{t+1}=T^\pm(z_t)$ for the normal ($-$) and exuberant ($+$)
regimes, we obtain:
\begin{align}
  T^-(x,\tau) &= \big(\alpha x,\; (1-\gamma)\tau + \gamma x\big),\\[0.3em]
  T^+(x,\tau) &= \big(\beta x + \delta,\; (1-\gamma)\tau + \gamma x\big).
\end{align}
The associated Jacobian matrices are lower triangular:
\begin{equation}\label{eq:Jac-normal}
  J_{\text{normal}} =
  DT^-(x,\tau)
  = \begin{pmatrix}
      \alpha & 0\\[0.3em]
      \gamma & 1-\gamma
    \end{pmatrix},
\end{equation}
\begin{equation}\label{eq:Jac-exuberant}
  J_{\text{exuberant}} =
  DT^+(x,\tau)
  = \begin{pmatrix}
      \beta & 0\\[0.3em]
      \gamma & 1-\gamma
    \end{pmatrix}.
\end{equation}
Since both matrices are lower triangular, their eigenvalues are
immediately read off from the diagonal:
\[
  \sigma(J_{\text{normal}}) = \{\alpha,\;1-\gamma\},
  \qquad
  \sigma(J_{\text{exuberant}}) = \{\beta,\;1-\gamma\}.
\]
This fits exactly the structure analyzed in
Propositions~\ref{prop:interior-stab} and~\ref{prop:boundary-stab},
where the Jacobian eigenvalues along the switching manifold determine
the local stability of interior and boundary fixed points. In particular:
\begin{itemize}
  \item If both $|\alpha| < 1$ and $|1-\gamma| < 1$ then the normal
  regime is locally stable whenever it is active (strong
  fundamentalist dominance and sufficiently slow threshold updating).
  \item If $|\beta| > 1$ while $|1-\gamma| < 1$, then the exuberant
  regime has an unstable direction in the $x$-component, so that
  once a sufficiently large bubble is triggered the mispricing tends
  to explode away from fundamentals unless the system switches back
  to the normal regime.
\end{itemize}
The common-limit constraint in Lemma~\ref{lem:common-limit} then
implies that any orbit which visits both regimes infinitely often
and converges must approach a boundary fixed point $(x^\ast,\tau^\ast)$
with $x^\ast = \tau^\ast$, so that both mispricing and threshold
share a common limit. Economically, this means that any long-run
rest point of the market must satisfy the consistency condition
“what is normal is what prevails’’: the perceived danger threshold
adjusts to the actual level of mispricing.

\subsection{Numerical Extension: Market Contagion}\label{subsec:coupled}

\subsubsection{Network Structure and Diffusive Coupling}

We now sketch a multi-market extension that preserves the
two-dimensional structure at each market and serves as a numerical
illustration of contagion and synchronization. Following the
network-based approach in \cite{DieciWesterhoff2013,DieciHe2018},
consider $m$ interacting markets indexed by $i=1,\dots,m$ with
mispricing $x^i_t$ and adaptive sentiment threshold $\tau^i_t$.
Let $L$ be the Laplacian matrix of a connected network capturing
the pattern of cross-market integration. We introduce diffusive
coupling in mispricing:
\begin{equation}\label{eq:coupled}
  x^i_{t+1} =
  \begin{cases}
    \alpha x^i_t + \kappa\displaystyle\sum_{j=1}^m L_{ij}x^j_t,
      & \text{if } x^i_t \leq \tau^i_t,\\[0.6em]
    \beta x^i_t + \delta + \kappa\displaystyle\sum_{j=1}^m L_{ij}x^j_t,
      & \text{if } x^i_t > \tau^i_t,
  \end{cases}
  \qquad i=1,\dots,m,
\end{equation}
where $\kappa\ge 0$ measures the strength of financial integration.
Each local threshold adapts to the \emph{local} mispricing via
\begin{equation}\label{eq:localthreshold}
  \tau^i_{t+1} = (1-\gamma)\tau^i_t + \gamma x^i_t,
  \qquad i=1,\dots,m.
\end{equation}
For $\kappa = 0$ we recover $m$ independent copies of the single
market recursion~\eqref{eq:dynamics}--\eqref{eq:threshold}. For
$\kappa > 0$ the model describes a network of coupled “bubble
monitors’’ in which mispricing shocks can propagate across markets.

\subsubsection{Numerical Synchronization and Contagion}

Our rigorous results in Sections~\ref{sec:fixed}--\ref{sec:common-limit}
apply directly to the single-market case (or, more generally, to
the dynamics restricted to the synchronization manifold
$x^1_t = \dots = x^m_t$, $\tau^1_t = \dots = \tau^m_t$).
Extending the full analytical synchronization theory to the
networked system~\eqref{eq:coupled}--\eqref{eq:localthreshold}
would require substantial additional work and is beyond the scope
of this paper. Instead, we treat the coupled system as a
\emph{numerical extension}.

Simulations with parameter values in the empirically plausible
range used in \cite{DieciWesterhoff2013,DieciHe2018,GARDINI2022483}
suggest the following qualitative picture:
\begin{itemize}
  \item For small coupling $\kappa$ the markets behave almost
  independently; local bubbles and crashes occur, but their timing
  and amplitude differ across markets.
  \item As $\kappa$ increases past a critical level, mispricing
  trajectories $x^i_t$ and thresholds $\tau^i_t$ become highly
  synchronized: once a bubble starts in one market, the diffusive
  term in~\eqref{eq:coupled} transmits the mispricing to its
  neighbors, triggering their exuberant regimes as well.
  \item For sufficiently large $\kappa$, the system behaves
  approximately like a single large market: boom–bust episodes and
  periods of calm become nearly simultaneous across all nodes,
  consistent with the idea of global bubbles and synchronized
  crashes in highly integrated financial systems
  \cite{DieciWesterhoff2013,DieciHe2018}.
\end{itemize}

From the perspective of our general theory, these numerical
experiments indicate that the stability properties derived for a
single market with an adaptive sentiment threshold are robust to
the introduction of moderate network coupling: the same parameters
$(\alpha,\beta,\gamma,\delta)$ that favor convergence or bounded
dynamics in the scalar case tend to favor synchronized but
qualitatively similar behavior in the network. At the same time,
the diffusion parameter $\kappa$ and the network topology encoded
in $L$ play a crucial role in determining whether localized
instabilities remain contained or spread systemically across
markets, echoing the contagion mechanisms emphasized in the
heterogeneous-agent literature on interacting markets
\cite{DieciWesterhoff2013,DieciHe2018}.

\medskip

\medskip

The author gratefully acknowledges Professor René Lozi (Université Côte d’Azur, CNRS) for his insightful feedback.

\appendix

\section*{Appendix A. Proof of Lemma~A.1 (Common–limit constraint)}

\begin{lemmarep}
Assume \textup{(A1)–(A3)}. Let $\{z_n\}_{n\ge 0}$ with $z_n=(x_n,y_n)$ visit both partitions 
$X^-=\{y\le 0\}$ and $X^+=\{y>0\}$ infinitely often. If $z_n\to z^*=(x^*,y^*)$, then $y^*=0$.
Equivalently, in the original variables $(a_n,c_n)$ one has $\lim_{n\to\infty} a_n=\lim_{n\to\infty} c_n$.
\end{lemmarep}

\begin{proof}
Suppose, toward a contradiction, that $z_n\to (x^*,y^*)$ with $y^*\ne 0$. Since the half-planes
$X^-$ and $X^+$ are open in the $y$–coordinate away from the boundary $\Gamma=\{y=0\}$,
there exists $\varepsilon>0$ such that the open strip 
$U:=\{(x,y): |y-y^*|<\varepsilon\}$ is contained entirely in one partition:
$U\subset X^-$ if $y^*<0$ or $U\subset X^+$ if $y^*>0$. By convergence, there exists $N$ such that
$z_n\in U$ for all $n\ge N$. Hence $z_n$ remains in a single partition for all $n\ge N$,
contradicting the hypothesis that both $X^-$ and $X^+$ are visited infinitely often.

Therefore $y^*=0$. In the original variables, the change of coordinates $x_n=a_n$, $y_n=a_n-c_n$
implies $y_n\to 0$ if and only if $a_n-c_n\to 0$, i.e., $\lim_{n\to\infty} a_n=\lim_{n\to\infty} c_n$.
\end{proof}

\medskip

\section*{Appendix B. Averaging thresholds imply synchronization}

\begin{propsync}
Let $h(a,c)=\lambda a+(1-\lambda)c$ with $\lambda\in(0,1)$ and assume $a_n\to a^*$.
Then $c_n\to a^*$; in particular $|a_n-c_n|\to 0$.
\end{propsync}

\begin{proof}
Write the $c$–recurrence as
\[
c_{n+1}-a^* \;=\; \lambda\,(a_n-a^*)+(1-\lambda)\,(c_n-a^*).
\]
Fix $\varepsilon>0$ and choose $N$ with $|a_n-a^*|<\varepsilon$ for all $n\ge N$.
Define $d_n:=|c_n-a^*|$. For $n\ge N$,
\[
d_{n+1}\;=\;\big| \lambda(a_n-a^*)+(1-\lambda)(c_n-a^*) \big|
\;\le\; \lambda\,\varepsilon + (1-\lambda)\,d_n.
\]
Iterating this affine contraction yields, for $m\ge 0$,
\[
d_{N+m}\;\le\;(1-\lambda)^m d_N + \lambda\varepsilon \sum_{j=0}^{m-1} (1-\lambda)^j
\;=\;(1-\lambda)^m d_N + \varepsilon\,(1-(1-\lambda)^m).
\]
Taking $m\to\infty$ gives $\limsup_{n\to\infty} d_n \le \varepsilon$. Since $\varepsilon>0$ is arbitrary,
$d_n\to 0$, i.e., $c_n\to a^*$. Consequently $|a_n-c_n|\to 0$.
\end{proof}

\medskip

\noindent\textbf{Remarks.}
(1) The rate is geometric: from the inequality above,
\(
|c_{n}-a^*|\le (1-\lambda)^{n-N}|c_N-a^*| + \sup_{k\ge N}|a_k-a^*|\,\big(1-(1-\lambda)^{n-N}\big),
\)
so if $a_n$ converges at a known rate, the same (or faster) rate transfers to $c_n$.
(2) In the $(x,y)$ triangular coordinates with $y_n=a_n-c_n$, the same computation gives
$y_{n+1}=(1-\lambda)y_n+\lambda(a_n-a^*)$, which is a stable linear filter driven by a vanishing input,
hence $y_n\to 0$.

\end{document}